\documentclass[12pt, twoside, leqno]{article}

\usepackage{amsmath,amsthm}
\usepackage{amssymb}
\usepackage[all]{xy}
\usepackage{enumerate}

\usepackage{graphicx}

\newtheorem{thm}{Theorem}[section]
\newtheorem{cor}[thm]{Corollary}
\newtheorem{lem}[thm]{Lemma}

\theoremstyle{definition}
\newtheorem{defin}[thm]{Definition}
\newtheorem{rem}[thm]{Remark}

\numberwithin{equation}{section}

\begin{document}


\baselineskip=17pt


\title{On the Capacity and Depth of Compact Surfaces}

\author{M. Abbasi\\
E-mail: mahboubeh.abbasi@gmail.com
\and 
B. Mashayekhy\\
Department of Pure Mathematics,\\
Center of Excellence in Analysis\\
 on Algebraic Structures\\
Ferdowsi University of Mashhad\\
P.O.Box 1159-91775, Mashhad, Iran\\
E-mail: bmashf@um.ac.ir}

\date{}

\maketitle


\renewcommand{\thefootnote}{}

\footnote{2010 \emph{Mathematics Subject Classification}: 55P55, 55P15, 55P20.}

\footnote{\emph{Key words and phrases}: Capacity, Compact surface, Homotopy domination, Homotopy type, Eilenberg-MacLane space, Polyhedron.}

\renewcommand{\thefootnote}{\arabic{footnote}}
\setcounter{footnote}{0}


\begin{abstract}
K. Borsuk in 1979, at the Topological Conference in Moscow, introduced the concept of capacity and depth of a compactum. In this paper we compute the capacity and depth of compact surfaces. We show that the capacity and depth of every compact orientable surface of genus $g\geq 0$ is equal to $g+2$. Also, we prove that the capacity and depth of a compact non-orientable surface of genus $g>0$ is $[\frac{g}{2}]+2$.
\end{abstract}

\section{Introduction}
K. Borsuk in \cite{So} introduced the concept of capacity and depth of a compactum (compact metric space) as follows:
 the capacity $C(A)$ of a compactum $A$ is the cardinality of the set of all shapes of compacta $X$ such that $\mathcal{S}
h(X) \leqslant \mathcal{S}h(A)$.
A system $\mathcal{S}h(X_{1}) < \mathcal{S}h(X_{2}) < \cdots < \mathcal{S}h(X_{k}) \leqslant \mathcal{S}h(A)$ is called a chain of length $k$ for compactum $A$, where $\mathcal{S}h(X_{i}) < \mathcal{S}h(X_{i+1})$ denotes $\mathcal{S}h(X_{i}) \leqslant \mathcal{S}h(X_{i+1})$ but not $\mathcal{S}h(X_{i+1}) \leqslant \mathcal{S}h(X_{i})$. The depth $D(A)$ of a compactum $A$ is the least upper bound of the lengths of all chains for $A$. It is clear that $D(A) \leq C(A)$ for each compactum $A$.

In the case of polyhedra,  the notions shape and shape domination in the above definitions can be replaced by the notions homotopy type and homotopy domination, respectively. Indeed, by some known
results in shape theory one can conclude that for any polyhedron $P$, there is a 1-1 functorial correspondence between the shapes of compacta shape dominated by $P$ and the homotopy types of CW-complexes (not necessarily finite) homotopy dominated by $P$ (in both pointed and unpointed cases) \cite{18}.

S. Mather in \cite{17} proved that every polyhedron dominates only countably many of different homotopy types (hence shapes). Since the capacity of a topological space is a homotopy invariant, i.e., if topological spaces $X$ and $Y$ have the same homotopy type, then $C(X)=C(Y)$, it is interesting to know which topological spaces have finite capacity and compute their exact capacity.  Borsuk in \cite{So} asked a question: `` Is it true that the capacity of every finite polyhedron is finite? ''.  D. Kolodziejczyk in \cite{16} gave a negative answer to this question. Also, in \cite{13} she proved that there exist polyhedra with infinite capacity and finite depth. Moreover, she investigated some conditions for polyhedra to have finite capacity (see \cite{18,14,15}). For instance, polyhedra with finite fundamental groups and polyhdera $P$ with abelian fundamental groups $\pi_1 (P)$ and finitely generated homology groups $H_i (\tilde{P})$, for $i\geq 2$, have finite capacity. Also, in \cite{19} she proved that polyhedra with virtually polycyclic fundamental group have finite depth.

Borsuk in \cite{So} mentioned that the capacity of $\bigvee_k \mathbb{S}^1$ and $\mathbb{S}^n$ are equal to $k+1$ and 2, respectively.  It is known that every compactum shape dominated by a polyhedron $P$ has the shape of some $P$-like compactum (see, for example, \cite{Kod}). Therefore, the capacities of some CW-complexes, such as $\mathbb{R}\mathbb{P}^{n}$, $\mathbb{C}\mathbb{P}^{n}$, and $\mathbb{T}^{n}$ ($n$-dimensional torus) may be obtained from shape-theoretical results concerning $P$-like compacta (see \cite{Hand,Eber,Kees,Kod}).

In this paper we compute capacity and depth of compact surfaces. In fact, we show that the capacity and depth of every compact orientable surface of genus $g\geq 0$ is equal to $g+2$. Also, we prove that the capacity and depth of a compact non-orientable surface of genus $g>0$ is $[\frac{g}{2}]+2$.

\section{Preliminaries}
In this paper every topological space is assumed to be connected. We expect that the reader is familiar with the basic notions and facts of shape theory (see \cite{B1} and \cite{Mar}) and retract theory (see \cite{ret}). We need the following results and definitions for the rest of the paper.

\begin{defin}\cite{18}.
A homomorphism $g:G \longrightarrow H$ of groups is an \emph{$r$-homomorphism} if there
exists a homomorphism $f :H\longrightarrow G$ such that $g\circ f = id_H$. In this case $H$ is called an $r$-image of $G$.
\end{defin}
\begin{defin}\cite{Wu}
 Let $G$ be a group with a subgroup $H$. Then $H$ is called a \emph{retract} of $G$ if there exists a homomorphism $r:G\longrightarrow H$ such that $r\circ i =id_H$, where $i:H\longrightarrow G$ is the inclusion homomorphism.
\end{defin}
Note that if a group $H$ is an $r$-image of $G$ with an $r$-homomorphism $g$ and the converse homomorphism $f$, then $H\cong f(H)$ and $f(H)$ is a retract of $G$.

\begin{defin}\cite{Wu}.
The group $G$ is called a surface group if $G\cong \pi_1 (\mathcal{S})$ for a closed (compact without boundary) surface $\mathcal{S}$ with $\chi (\mathcal{S})<0$, where $\chi (\mathcal{S})$ is the Euler characteristic of $\mathcal{S}$.
\end{defin}
Recall that the Euler Characteristic of connected sum of $n$ tori and  connected sum of $n$ real projective planes are $2-2n$ and $2-n$, respectively.
\begin{lem}\cite[Lemma 4.5]{Wu}.\label{3}
Let $G$ be a surface group. If $K$ is any proper retract of $G$, then $K$ is
a free group with $rank\; K\leq \frac{1}{2}\; rank\; G$, where $rank\; G$ is the minimal number of generators of $G$.
\end{lem}

\begin{rem}\cite{Use}.\label{use}
Let $G=H\ltimes _\varphi U$ be the semidirect product of $U$ and $H$ with respect to $\varphi : H\rightarrow {\rm Aut}(U)$. For a subgroup $\Gamma$ of $G$, put $U_{\Gamma}= \Gamma\cap U$ and $\Gamma_{H}=\lbrace h\in H\; | \; (u, h)\in \Gamma\ {\rm for}\ {\rm some}\ u\in U\rbrace$. Group monomorphisms $\rho:R \longrightarrow H$ and $\lambda:L \longrightarrow U$ are called a $\varphi^{H}_{U}-pair$ of groups if there exists a short exact sequence of groups
$$
1\rightarrow L \xrightarrow{\alpha} T \xrightarrow{\beta} R \rightarrow 1
$$
with a monomorphism $\mu:T\longrightarrow G$ such that the following diagram commutes
\begin{displaymath}
\xymatrix{
L \ar[r]^\alpha \ar[rd]_\lambda &
T \ar[d]_\mu \ar[r]^\beta &
R \ar[ld]^\rho\\
&H\ltimes_\varphi U.}
\end{displaymath}
We denote a $\varphi^{H}_{U}-pair$ of groups by $ \langle L^{\lambda}, R^{\rho}, \varphi^{H}_{U}\rangle$. Usenko proved that each  $\varphi^{H}_{U}-pair$ of groups determines a subgroup of the semidirect product $G = H\ltimes_\varphi U$. On the other hand, each subgroup $\Gamma \leqslant G$ determines some $\varphi^{H}_{U}-pair$ of groups $ \langle U_{\Gamma}^{i}, \Gamma_{H}^{j}, \varphi^{H}_{U}\rangle$, where $i: U_{\Gamma}\longrightarrow U$ and $j: \Gamma_{H}\longrightarrow H$ are embeddings.
\end{rem}

\section{The Capacity of Compact Surfaces}
In this section we compute the capacity of compact (orientable or non-orientable) surfaces.
By the well-known classification theorem for compact surfaces, any compact orientable surface is homeomorphic to a sphere, or to a connected sum of $n$ tori, and any compact non-orientable surface is homeomorphic to a connected sum of $n$ real projective planes (see, for example, \cite[Theorem 7.2]{Mas}).
Recall that the capacity of $\mathbb{S}^2$ is equal to 2.

\begin{thm}\label{7}
The capacity of connected sum of $n$ tori is equal to $n+2$.
\end{thm}
\begin{proof}
Suppose that $\mathcal{S}$ is a connected sum of $n$ tori. If $n=1$, then one could easily conclude that the capacity of the torus $\mathbb{T}^{2}$ is $3$. Let $n>1$ and put $G=\pi_1 (\mathcal{S})$. Suppose that a space $\mathcal{A}$ is homotopy dominated by $\mathcal{S}$. Then $\pi_1 (\mathcal{A})$ is isomorphic to a retract of $G$. We have the following two cases:

Case One. $\pi_1 (\mathcal{A})$ is isomorphic to a proper retract of $G$. Since $n >1$,  $\chi (\mathcal{S})<0$, and hence $G$ is a surface group. Because $G$ is a finitely presented group with $2n$ generators and one relation (see \cite[p.\;99]{Mas}), $\pi_1 (\mathcal{A})$ is a free group with $rank\; \pi_1 (\mathcal{A})\leq \frac{1}{2}\; rank\; G=n$ by
using Lemma \ref{3}.
Now since $\mathcal{S}$ (and hence $\mathcal{A}$) is an Eilenberg-MacLane space (by \cite[Example 1B.2]{3}, every orientable or non-orientable compact surface of genus $g>1$, is an Eilenberg-MacLane space $K(G,1)$), $\mathcal{A}$ has the form of one of the following spaces:
\[
*, \; \mathbb{S}^1 , \; \mathbb{S}^1 \vee \mathbb{S}^1 , \cdots , \; \underbrace{\mathbb{S}^1 \vee \cdots \vee \mathbb{S}^1}_{n-folds}.
\]
On the other hand, $ \underbrace{\mathbb{S}^1 \vee \cdots \vee \mathbb{S}^1}_{n-folds}$ is a retract of $\mathcal{S}$, and so each of the spaces above is a retract of $\mathcal{S}$.

Case Two. $\pi_1 (\mathcal{A})\cong G$. The spaces $\mathcal{A}$ and $\mathcal{S}$ have the same homotopy type by the fact that the fundamental group of $\mathcal{S}$ determines the homotopy type of $\mathcal{S}$.
\end{proof}

\begin{thm}\label{Kl}
The capacity of the Klein bottle is equal to 3.
\end{thm}
\begin{proof}
Suppose that $\mathcal{K}$ is the Klein bottle. We know that $\pi_1 (\mathcal{K})$ has a presentation as $\langle x,y\; | \; yxy^{-1}=x^{-1} \rangle \cong \mathbb{Z}\ltimes\mathbb{Z}$ (see \cite[p.\;79]{Mas}). Suppose that $\Gamma \leqslant \mathbb{Z} \ltimes \mathbb{Z}$. By the argument of Remark \ref{use}, $\Gamma \cong \Gamma_{H} \ltimes U_{\Gamma}$ since $U_{\Gamma}$ (and also $\Gamma_{H}$) is trivial or $\mathbb{Z}$, up to isomorphism. Hence isomorphism classes of subgroups of $\mathbb{Z} \ltimes \mathbb{Z}$ are trivial, $\mathbb{Z}$, $\mathbb{Z}\times\mathbb{Z}$ and $\mathbb{Z}\ltimes \mathbb{Z}$. Since $\pi_1 (\mathcal{K})$ determines the homotopy type of $\mathcal{K}$, the spaces $*$, $\mathbb{S}^1$, and $\mathcal{K}$ (which are $K(1,1)$, $K(\mathbb{Z},1)$, and $K(\mathbb{Z}\ltimes \mathbb{Z},1)$, respectively) are homotopy dominated by $\mathcal{K}$. But $K(\mathbb{Z}\times \mathbb{Z},1)$ is not homotopy dominated by $\mathcal{K}$ because $H_1 (K(\mathbb{Z}\times \mathbb{Z},1) )\cong \mathbb{Z}\times \mathbb{Z}$ is not a retract of $H_1 (\mathcal{K})\cong \mathbb{Z}\times \mathbb{Z}_2$. Hence the proof is complete.
\end{proof}
\begin{lem}\label{conp}
Let $\mathcal{S}$ be a connected sum of $n$ real projective planes, where $n>2$. Then $\underbrace{\mathbb{S}^1 \vee \cdots \vee \mathbb{S}^1}_{[\frac{n}{2}]-folds}$ is a retract of $\mathcal{S}$.
\end{lem}
\begin{proof}
We consider the following two cases:

Case One: $n=2m+1$. Then $\mathcal{S}$ is homeomorphic to connected sum of a real projective plane and  $m$ tori (see \cite[p. 101]{Mas}). 
First, we show that $\underbrace{\mathbb{S}^1 \vee \cdots \vee \mathbb{S}^1}_{m-folds}$ is a retract of $\underbrace{\mathbb{T}^2 \# \cdots \# \mathbb{T}^2}_{m-folds}$. Although one can write equations for the retraction $r_1:\underbrace{\mathbb{T}^2 \# \cdots \# \mathbb{T}^2}_{m-folds}\longrightarrow \underbrace{\mathbb{S}^1 \vee \cdots \vee \mathbb{S}^1}_{m-folds}$, it is simpler to indicate it in picture, as we have done in Figure  \ref{fig1}. Let $\mathcal{M}$ be the wedge sum of $m$ tori. One maps $\underbrace{\mathbb{T}^2 \# \cdots \# \mathbb{T}^2}_{m-folds}$ onto $\mathcal{M}$ by a map that collapses each the dotted circle to a point but is otherwise one-to-one. Therefore, it defines a homeomorphism $h$ of the $\underbrace{\mathbb{S}^1 \vee \cdots \vee \mathbb{S}^1}_{m-folds} \subseteq \underbrace{\mathbb{T}^2 \# \cdots \# \mathbb{T}^2}_{m-folds}$ with the 
$\underbrace{\mathbb{S}^1 \vee \cdots \vee \mathbb{S}^1}_{m-folds} \subseteq \mathcal{M}$. Then one retracts $\mathcal{M}$ onto its $\underbrace{\mathbb{S}^1 \vee \cdots \vee \mathbb{S}^1}_{m-folds}$ by mapping each cross-sectional circle to the point where it intersects $\underbrace{\mathbb{S}^1 \vee \cdots \vee \mathbb{S}^1}_{m-folds}$.
Finally, one maps the $\underbrace{\mathbb{S}^1 \vee \cdots \vee \mathbb{S}^1}_{m-folds}$ in $\mathcal{M}$ back onto the $\underbrace{\mathbb{S}^1 \vee \cdots \vee \mathbb{S}^1}_{m-folds}$ in $\underbrace{\mathbb{T}^2 \# \cdots \# \mathbb{T}^2}_{m-folds}$ by the map $h^{-1}$ (compare \cite[p.\;374]{Munk}).
\begin{figure}
\centering
\includegraphics[scale=0.4]{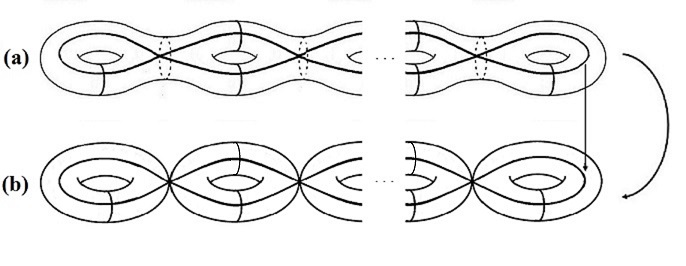}
\caption{\label{fig1}(a) A connected sum of m tori; (b) The space $\mathcal{M}$}
\end{figure}
Now assume that  $C_1$ is the common circle of connected sum of $\underbrace{\mathbb{T}^2 \# \cdots \# \mathbb{T}^2}_{m-folds}$ and $\mathbb{RP}^2$. Without loss of generality, suppose that $ C_{1} \bigcap \underbrace{\mathbb{S}^1 \vee \cdots \vee \mathbb{S}^1}_{m-folds} = \lbrace x_{0} \rbrace$. 

Let $q:\underbrace{\mathbb{T}^2 \# \cdots \# \mathbb{T}^2}_{m-folds}\# \mathbb{RP}^2\longrightarrow \underbrace{\mathbb{T}^2 \# \cdots \# \mathbb{T}^2}_{m-folds}\vee \mathbb{RP}^2$ be the quotient map such that collapses $C_{1}$ to the point $\lbrace x_{0} \rbrace$.

Also, the map $r_2:\underbrace{\mathbb{T}^2 \# \cdots \# \mathbb{T}^2}_{m-folds}\vee \mathbb{RP}^2 \longrightarrow \underbrace{\mathbb{T}^2 \# \cdots \# \mathbb{T}^2}_{m-folds}$ defined by $r_2(x,y)=x$ is a retraction (here the space $\underbrace{\mathbb{T}^2 \# \cdots \# \mathbb{T}^2}_{m-folds}\vee \mathbb{RP}^2$ is considered as subspace $\underbrace{\mathbb{T}^2 \# \cdots \# \mathbb{T}^2}_{m-folds} \times  \lbrace y_{0}\rbrace  \bigcup \lbrace x_{0}\rbrace \times \mathbb{RP}^2$ of space $\underbrace{\mathbb{T}^2 \# \cdots \# \mathbb{T}^2}_{m-folds} \times \mathbb{RP}^2$, where $x_0 \in \underbrace{\mathbb{T}^2 \# \cdots \# \mathbb{T}^2}_{m-folds}$ and $y_0 \in \mathbb{RP}^2$). Hence the continuous map 
\[ 
r_1\circ r_2 \circ q:\underbrace{\mathbb{T}^2 \# \cdots \# \mathbb{T}^2}_{m-folds}\# \mathbb{RP}^2\longrightarrow \underbrace{\mathbb{S}^1 \vee \cdots \vee \mathbb{S}^1}_{m-folds}
\]
is  a retraction.

Case Two:  $n=2m$. Then $\mathcal{S}$ is homeomorphic to connected sum of a Klein bottle and  $m-1$ tori (see \cite[p. 101]{Mas}). Assume that $q:\underbrace{\mathbb{T}^2 \# \cdots \# \mathbb{T}^2}_{(m-1)-folds}\# \mathcal{K}\longrightarrow \underbrace{\mathbb{T}^2 \# \cdots \# \mathbb{T}^2}_{(m-1)-folds}\vee  \mathcal{K}$ is the quotient map (which is defined by a similar definition of map $q$ in Case One). Observe that $q(x)=x$, for all $ x \in \Big(\underbrace{\mathbb{S}^1 \vee \cdots \vee \mathbb{S}^1}_{(m-1)-folds}\Big) \vee \mathbb{S}^1$, where $\underbrace{\mathbb{S}^1 \vee \cdots \vee \mathbb{S}^1}_{(m-1)-folds} \subseteq \underbrace{\mathbb{T}^2 \# \cdots \# \mathbb{T}^2}_{(m-1)-folds}$ and $\mathbb{S}^1 \subseteq \mathcal{K}$. Also, there exist retractions $r_1 :\underbrace{\mathbb{T}^2 \# \cdots \# \mathbb{T}^2}_{(m-1)-folds}\longrightarrow \underbrace{\mathbb{S}^1 \vee \cdots \vee \mathbb{S}^1}_{(m-1)-folds}$ and $r_2 :\mathcal{K}\longrightarrow \mathbb{S}^1$ (by the proof of Theorem \ref{Kl}).

Now we define map $r:\underbrace{\mathbb{T}^2 \# \cdots \# \mathbb{T}^2}_{(m-1)-folds}\vee \mathcal{K} \longrightarrow \underbrace{\mathbb{S}^1 \vee \cdots \vee \mathbb{S}^1}_{m-folds}$ as follows:
\[
r(x)=
\begin{cases}
r_1 (x), & q(x)\in \underbrace{\mathbb{T}^2 \# \cdots \# \mathbb{T}^2}_{(m-1)-folds}\\
r_2 (x), & 	q(x)\in \mathcal{K}.
\end{cases}
\]
 Let $C_2$ be the common circle of connected sum of $\underbrace{\mathbb{T}^2 \# \cdots \# \mathbb{T}^2}_{(m-1)-folds}$ and $\mathcal{K}$. For all $x\in C_2$, we have $r_i (q(x))=r_i (*)=*$, where $i=1,2$ and $*$ is the common point of the wedge sum. One can see that $r$ is continuous. It is enough to show that $r$ is a retraction. If $x\in \underbrace{\mathbb{S}^1 \vee \cdots \vee \mathbb{S}^1}_{(m-1)-folds}$ ($\subseteq \underbrace{\mathbb{T}^2 \# \cdots \# \mathbb{T}^2}_{(m-1)-folds}$), then $r(x)=r_1 (x)=x$, and if $x\in \mathbb{S}^1$ ($\subseteq \mathcal{K}$), then $r(x)=r_2 (x)=x$. Therefore, $\underbrace{\mathbb{S}^1 \vee \cdots \vee \mathbb{S}^1}_{m-folds}$ is a retract of $\underbrace{\mathbb{T}^2 \# \cdots \# \mathbb{T}^2}_{(m-1)-folds}\vee \mathcal{K}$. Hence the continuous map $r\circ q: \underbrace{\mathbb{T}^2 \# \cdots \# \mathbb{T}^2}_{(m-1)-folds}\# \mathcal{K}\longrightarrow \underbrace{\mathbb{S}^1 \vee \cdots \vee \mathbb{S}^1}_{m-folds}$ satisfies $r \circ q(x)=x$, for all $x \in \underbrace{\mathbb{S}^1 \vee \cdots \vee \mathbb{S}^1}_{m-folds}$. 
 
Accordingly, $\underbrace{\mathbb{S}^1 \vee \cdots \vee \mathbb{S}^1}_{m-folds}$ is a retract of $\underbrace{\mathbb{RP}^2 \# \cdots \# \mathbb{RP}^2}_{n-folds}$, where $m=[\frac{n}{2}]$.
\end{proof}
\begin{thm}\label{13}
The capacity of connected sum of $n$ real projective planes is equal to $[\frac{n}{2}]+2$.
\end{thm}
\begin{proof}
Suppose that $\mathcal{S}$ is the connected sum of $n$ real projective planes and $\mathcal{A}$ is homotopy dominated by $\mathcal{S}$.

If $n=1$, then the fact that $C(\mathbb{RP}^{2}) = 2$ is easily concluded (see \cite{Kod} or \cite{Hand}). If $n=2$, then $\mathcal{S}$ is a Klein bottle (see \cite[Example 4.3]{Mas}), and so by Theorem \ref{Kl} the capacity of $\mathcal{S}$ is 3. Now suppose that $n>2$ and put $G=\pi_1 (\mathcal{S})$. 
We have the following two cases:

Case One:  $n=2m+1$ where $m \geq 1$. Then $\mathcal{S}$ is homeomorphic to connected sum of a real projective plane and  $m$ tori. Since $n > 2$, $\chi (\mathcal{S})<0$, and so $G$ is a surface group.
It is obvious that $\pi_1 (\mathcal{A})$ is isomorphic to a retract of $G$. If $\pi_1 (\mathcal{A})$ is isomorphic to a proper retract of $G$, then by Lemma \ref{3} and the fact that $G$ is a finitely presented group with $n$ generators and one relation (see \cite[p.\;100]{Mas}), $\pi_1 (\mathcal{A})$ is a free group with $rank\; \pi_1 (\mathcal{A})\leq \frac{1}{2}\; rank\; G=\frac{n}{2}$.  

Since $\mathcal{S}$ (and so $\mathcal{A}$) is an Eilenberg-MacLane space and $\pi_1 (\mathcal{A})$ is a free group of rank $0\leq t\leq [\frac{n}{2}]=m$, the space $\mathcal{A}$ has the form of one of the following spaces:
\[
*, \: \mathbb{S}^1 , \: \mathbb{S}^1 \vee \mathbb{S}^1 , \cdots , \: \underbrace{\mathbb{S}^1 \vee \cdots \vee \mathbb{S}^1}_{m-folds}.
\]
On the other hand, by Lemma \ref{conp}, $\underbrace{\mathbb{S}^1 \vee \cdots \vee \mathbb{S}^1}_{m-folds}$ is a retract of $\mathcal{S}$, and so each of the spaces above can be a retract of $\mathcal{S}$.

Now if $\pi_1 (\mathcal{A})\cong G$, then $\mathcal{A}$ and $\mathcal{S}$ have the same homotopy type since $\mathcal{S}$ is an Eilenberg-MacLane space ,and so is $\mathcal{A}$.
Hence the capacity of $\mathcal{S}$ is equal to $m+2=[\frac{n}{2}]+2$.

Case Two:  $n=2m$ where $m > 1$. Then $\mathcal{S}$ is homeomorphic to connected sum of a Klein bottle and  $m-1$ tori.  By a similar argument to the case one, $\mathcal{A}$ has  the form of one of the following spaces:
\[
*, \: \mathbb{S}^1 , \: \mathbb{S}^1 \vee \mathbb{S}^1 , \cdots , \: \underbrace{\mathbb{S}^1 \vee \cdots \vee \mathbb{S}^1}_{m-folds} , \mathcal{S}.
\]
By Lemma \ref{conp}, $\underbrace{\mathbb{S}^1 \vee \cdots \vee \mathbb{S}^1}_{m-folds}$ is a retract of $\mathcal{S}$. Thus each of the spaces above can be a retract of $\mathcal{S}$.
Hence the capacity of $\mathcal{S}$ is equal to $m+2=[\frac{n}{2}]+2$.
\end{proof}

Using Theorems \ref{7} and \ref{13} the following result can be easily concluded. 
\begin{cor}
The depth of orientable compact surfaces of genus $g \geq 0$ is $g+2$. Moreover, the depth of non-orientable compact surfaces of genus $g \geq 1$ is $[\frac{g}{2}]+2$.
\end{cor}


\end{document}